\documentclass[10pt,twoside,final]{camc}
\usepackage{enumerate}
\DeclareSymbolFont{igrfm}{OML}{cmm}{m}{it}
\def\camcsym#1#2{\DeclareMathSymbol{#1}{\mathord}{igrfm}{"#2}}
\camcsym{\alpha}{0B}\camcsym{\beta}{0C}\camcsym{\gamma}{0D}\camcsym{\delta}{0E}
\camcsym{\epsilon}{0F}\camcsym{\zeta}{10}\camcsym{\eta}{11}\camcsym{\theta}{12}
\camcsym{\iota}{13}\camcsym{\kappa}{14}\camcsym{\lambda}{15}\camcsym{\mu}{16}
\camcsym{\nu}{17}\camcsym{\xi}{18}\camcsym{\pi}{19}\camcsym{\rho}{1A}
\camcsym{\sigma}{1B}\camcsym{\tau}{1C}\camcsym{\upsilon}{1D}\camcsym{\phi}{1E}
\camcsym{\chi}{1F}\camcsym{\psi}{20}\camcsym{\omega}{21}\camcsym{\varepsilon}{22}
\camcsym{\vartheta}{23}\camcsym{\varpi}{24}\let\varrho\rho\let\varsigma\sigma
\camcsym{\varphi}{27}\camcsym{\Gamma}{00}\camcsym{\Delta}{01}\camcsym{\Theta}{02}
\camcsym{\Lambda}{03}\camcsym{\Xi}{04}\camcsym{\Pi}{05}\camcsym{\Sigma}{06}
\camcsym{\Upsilon}{07}\camcsym{\Phi}{08}\camcsym{\Psi}{09}\camcsym{\Omega}{0A}
\let\le\leqslant\let\leq\leqslant\let\ge\geqslant\let\geq\geqslant
\let\ldots\cdots\let\emptyset\varnothing
\newcommand{\edmt}[3][\footnotesize]{
\par\noindent{#1\sffamily\bfseries #2~}{#1\rmfamily #3}\vskip.6\baselineskip}
\def\acknowl#1{\edmt{Acknowledgements}{#1}}
\def\funding#1{\edmt{Funding}{#1}}

\def\complia#1{\edmt[]{Compliance with Ethical Standards}{}\edmt{Conflict of Interest}{#1}}
\def\ethical#1{\edmt{Ethical Approval}{#1}}
\AtBeginDocument{%
\nocite{label}
\newcommand{\dx}{{\mathrm d}x}
\newcommand{\ds}{{\mathrm d}s}
\newcommand{\vbf}[1]{\textbf{\itshape#1}}
}

\usepackage{amssymb,amsmath,amsthm}
\usepackage{lineno}
\usepackage{xcolor}
\usepackage[colorlinks=true]{hyperref}

\renewcommand{\r}{\textsc r} 
\renewcommand{\j}{\textsc j} 
\renewcommand{\c}{\textsc c} 
\renewcommand{\o}{\textsc o} 
\renewcommand{\d}{\textsc d} 

\newcommand{\radius}{\mathcal R}
\newcommand{\KNOW}{K_\mathcal N}
\newcommand{\Dt}{\Delta t}
\newcommand{\vmax}{v_{\text{max}}}
\newcommand{\Dmem}{\delta m}
\newcommand{\Dcom}{\delta c}
\newcommand{\nextroad}{\textsc {nextroad}} 
\newcommand{\nextcar}{\textsc {nextcar}}   
\newcommand{\Nj}{N_\textsc j}
\newcommand{\Nr}{N_\textsc r}
\newcommand{\Nc}{N_\textsc c}
\newcommand{\TT}{w}
\newcommand{\TTT}{\textsc{ttt}}
\newcommand{\TTc}{\textsc{tt}_\c}
\newcommand{\Tfin}{T_{fin}}

\newcommand{\REV}[1]{\textcolor{black}{#1}}

\graphicspath{{all_figs/}}

\begin{document}

\title{A Microscopic Traffic Flow Model on Network with Destination-Aware V2V Communications and Rational Decision-Making}
\author[addressref={aff1}, corref, 
email={e.cristiani@iac.cnr.it}]{Emiliano Cristiani}%
\author[addressref={aff1,aff2},]{Francesca L. Ignoto}%
\address[id=aff1]{Istituto per le Applicazioni del Calcolo, Consiglio Nazionale delle Ricerche, Rome, Italy}%
\address[id=aff2]{Dipartimento di Scienze di Base e Applicate per l’Ingegneria, Sapienza Università di Roma, Rome, Italy}%
\maketitle
\abstract{%
In this paper we carry out a computational study of a novel microscopic follow-the-leader model for traffic flow on road networks. 
We assume that each driver has its own origin and destination, and wants to complete its journey in minimal time. We also assume that each driver is able to take rational decisions at junctions and can change route while moving depending on the traffic conditions. 
The main novelty of the model is that vehicles can automatically and anonymously share information about their position, destination, and planned path when they are close to each other within a certain distance. 
The pieces of information acquired during the journey are used to optimize the route itself.
In the limit case of \REV{an infinite} communication range, we recover the classical Reactive User Equilibrium and Dynamic User Equilibrium.
}
\keywords{Traffic flow modeling $\cdot$ V2V communications $\cdot$ differential games $\cdot$ optimal control problems}
\classification{
91A80 
$\cdot$ 
34H05 
$\cdot$
76A30 
$\cdot$
34B45  	
$\cdot$
90B18  	
$\cdot$
90C39  	
}

%
%
%
%
%
%



\section{Introduction}\label{sec:intro}  

\paragraph{Context and relevant literature}
This paper is devoted to the study of a \textit{microscopic} (agent-based) differential model for traffic flow on networks, where vehicles are able to communicate with each other within a certain distance. 
Note that we do not consider other than inter-vehicle communications, i.e.\ we assume that vehicles cannot send/receive information to/from the road infrastructure or any centralized server. 
We also assume that each vehicle has an assigned \textit{origin-destination} (OD) pair, and all vehicles start moving from their origin \textit{at the same time}. 
Moreover, each vehicle is equipped with an \REV{On-board Processing Unit} (OPU), which stores information about the geometry of the network, sends/receives information about position, destination, planned path, and is able to suggest rational decisions to the drivers regarding the minimum-time route towards destination. 
Since \textit{all} vehicles are assumed to have such capabilities, we fall into a \textit{noncooperative game-theoretic} approach for the overall traffic dynamics.

The network is assumed to be composed by a set of unidirectional \textit{roads} and \textit{junctions}. 
Note that, mathematically speaking, a network resembles a \textit{graph}, where \textit{arcs} coincide with the roads and \textit{nodes} coincide with the junctions, but here roads have a real physical extension, with their own space coordinates and are not a mere logical connections between nodes. 

Finally, we assume that roads are 1-dimensional (no lanes) and junctions are 0-dimensional, with no given priorities or traffic lights. 
Vehicles are 0-dimensional as well, even if a positive length will be actually assigned for theoretical purposes. 

\medskip

Mathematical models for traffic flow on a single road have a very long story dating back to the 1950s, with a huge amount of research papers published on the topic. 
\REV{It is useful to note that} the mathematical and the engineering literature often proceeded in parallel without meeting, and arriving at the same conclusions from different points of view. 
Starting from the 1990s, the extensions to road networks became frequent, and now the literature covers extensively all possible approaches, including \textit{microscopic} models, either differential (based on ordinary differential equations) and nondifferential (e.g., cellular automata), and \textit{macroscopic} models, typically based on partial differential equations.

In this paper we are mainly interested in game-theoretic aspects of traffic \cite{ahmad2023IEEEA}, which are, in turn, related to the (rational) behavior of the drivers, especially in the context of connected and autonomous vehicles (CAVs).
It is possible -- and useful -- ordering possible dynamic traffic assignments (DTAs) on the basis of their degree of rationality:
\begin{itemize}
	\item \textbf{Basic behavior (BB)} 
	Vehicles follow the shortest route connecting their OD pair. The shortest route is computed once and for all, assuming empty roads and static travel times. While moving, vehicles rule their velocity/acceleration on the basis of traffic regulations and local congestion level.
	\item \textbf{Reactive user equilibrium (RUE)} Drivers know in real-time the congestion level on the whole network, obtained either by the knowledge of the single positions of all cars (microscale) or by the probability density function defined on all roads (macroscale). The congestion level is used to update the travel times of each road and continuously recompute the fastest route to destination. These are the same modeling assumptions of the Hughes model in the context of pedestrian dynamics \cite{huang2009TRB, hughes2002TRB}.
	\item \textbf{Dynamic user equilibrium (DUE)} Each driver knows other vehicles' OD pairs, can predict the dynamics of all other vehicles, and is able to compute the fastest route on the basis of the \textit{current} and \textit{future} traffic conditions. Since this ``predict \& optimize'' action is equally and independently taken by all drivers, we face a noncooperative game-theoretic problem. If a solution exists, it is a Nash-like (or Wardropian, in the traffic context) equilibrium, which corresponds to a scenario in which all drivers have chosen a route, and no driver can lower its time-to-destination unilaterally switching to a different route \cite{cristiani2015NHM, han2019NSE, morandi20244OR}. (An improvement is actually possible, but two or more drivers should agree to change their routes.)
	In a macroscopic setting, this scenario is also known as \textit{mean field game} in the mathematical literature. 
	Note also that optimization is not necessarily done with respect to the route choice (in other words, the behavior at junctions), but it can also be done with respect to the vehicle's velocity \cite{huang2020DCDS-B, huang2021TRC, mo2024TS}. 
    Past drivers' decisions can also be taken into account \cite{bagagiolo2021JDG}.
	\item \textbf{System optimum (SO)} In this case, it is assumed that all drivers accept to follow a route externally assigned to them by a centralized entity, usually the road manager. The primary goal here is to optimize some \textit{global} network performance (e.g., minimizing the total travel time, pollution emissions, or congestion levels). 
\end{itemize}
In the recent years, many researchers have proposed models which fall in between the aforementioned behaviors, often referred to models with partial knowledge or partial control, bounded rationality, multi-informed, or myopic (in space or time).
Let us mention, e.g., \cite{carrillo2016M3AS} for hybrid BB-RUE model, \cite{cristiani2023CMS, pan2012proc} for hybrid RUE-DUE models, and \cite{morandi20244OR} for hybrid DUE-SO model. 

Mixed scenarios, where some vehicles behave in a way and the remaining vehicles behave in another way, were also proposed; see, e.g., \cite{festa2023JOTA, samaranayake2018TS, siri2023CDC}. 
This is very common when one considers traffic scenarios composed by both human-driven vehicles (HDVs) and connected \& automated vehicles (CAVs).

Let us also mention that, in this framework, it is defined the so-called \textit{price of anarchy} (PoA), i.e.\ the ratio between the total travel cost under a user equilibrium and the total travel cost under the SO assignment \cite{belov2021JITS, morandi20244OR}. 

\medskip

The main topic of interest of this paper is that of \textit{vehicle-to-vehicle} (V2V) \textit{communications}, i.e.\ the possibility that vehicles share information whenever they get close enough to each other. 
As far as we know, traffic models including V2V technology were not studied so far in the mathematical literature, while they were extensively studied in the engineering literature since long time; see, e.g., \cite{enkelmann2003proc}.
In this respect, it could be useful for the reader to recall that V2V communications are also known as \textit{inter-vehicle communications} (IVC), or \textit{vehicular ad hoc networks} (VANET). 
The last term can also include vehicle-to-infrastructure (V2I) communications \cite{dutta2024JNSM}, which are wireless communications directed to a centralized server installed by the road manager.

V2V communications can be used for different goals, well reviewed in the recent papers \cite{elzorkany2020TOTJ, alsudani2023JMCOM, alabdouli2025SC, won2018chapter}. They can be divided in two macro-categories:
\begin{itemize}
	\item \textit{Local goals}. Communications are dedicated to avoid front-rear collisions (vehicles transmit the imminent intention of breaking to the following vehicles), coordinate  merging  and  lane-changing  activities, and to avoid collisions at junctions (vehicles reach an agreement on which should go first at the intersection); see, e.g., \cite{elzorkany2020TOTJ}. 
	Communications can be also used to support cooperative cruise control and platooning, aiming at, e.g., reducing fuel consumption or avoid instabilities in traffic flow.  
	Data packets usually contain vehicle's position, velocity, acceleration, and direction of travel.
	\item \textit{Global goals}. Communications are used to spread alert information across the whole road network, allowing drivers, even far away, to reroute for avoiding congested areas. 
	Data packets usually contain vehicle's position, degree of local congestion, travel times, presence of accidents with the timestamp of their occurrence, information about road closures; see, e.g., \cite{desouza2015proc, lakas2009proc, ohara2007proc, zardosht2017JTTE}. 
    In \cite{xu2011proc}, data packets also include planned route in order to predict the encounters between vehicles and then optimize communication events across the network.
\end{itemize}

It is also worth citing the papers \cite{bazzi2011proc, alobeidyeen2023VC, kim2017TRC} that investigate how fast the information (whatever it is) spreads across the road network hoping from one vehicle to another, and how this feature can be exploited for several purposes.

Finally, let us mention that V2V communications can be performed using a number of different technologies, each having its maximum transmission range and comes with a number of technical problems \cite{elzorkany2020TOTJ}. 
One for all, let us mention the \textit{broadcast storm problem}, which corresponds to an excessive broadcast packets overwhelming the network and can happen when many vehicles, close to each other, bounce the same information back and forth. 
Good reviews for this and other technical problems, including privacy issues and security issues (hacking), are \cite{dutta2024JNSM, elzorkany2020TOTJ, ullah2025IST}.

\medskip

Although this paper is confined to the microscopic scale, we will always keep an eye on possible macroscopic counterparts of the considered models. 
With this in mind, it is useful to recall here some important papers on \textit{many-particle limits}. 
The first-order microscopic differential model we use here was proposed in \cite{colombo2014RSMUP}. 
In the same paper, authors show that its many-particle limit coincides with the classical LWR model \cite{lighthill1955PRSLA,richards1956OR}, whose discrete version (when the Godunov numerical scheme is employed for discretization) is also known as CTM \cite{CTM} in the engineering literature. 
The same limit was further studied in the paper \cite{difrancesco2015ARMA}, and then in \cite{ancona2025DCDS}.
The model generalization to road networks, extensively used in this paper, was proposed in \cite{cristiani2016NHM}, together with some preliminary results about its many-particle limit. 
The authors show that the best candidate for the macroscale limit is the path-based model proposed in the twin papers \cite{bretti2014DCDS-S, briani2014NHM}, and further studied in \cite{dovetta2022SIMA}.

\paragraph{Main contribution}
In this paper we study the potential of sharing, in a V2V communications context, the \textit{destinations} of road users. 
In a real scenario, we obviously imagine an encrypted and anonymous communication protocol which only share data in aggregate form, to avoid privacy issues. 
Moreover, we assume that all computations are performed by the OPU, in such a way that human drivers, if present, do not have direct access to sensible information. 
That said, in this paper we will not address this kind of problems, keeping an ``applied mathematics'' point of view. 
Once the OPU has acquired desired destinations from the neighbor vehicles, it can use this information for the owner's advantage, e.g.\ for minimizing the time to destination. 
If only past and current information is used, we get a new kind of RUE with partial knowledge, while if we consider also future traffic conditions, predicted by the OPU itself, we get a new kind of DUE, again with partial knowledge. Both equilibria are modulated by the transmission range, which acts as a parameter, and they are expected to converge, respectively, to the standard RUE and to the standard DUE in the limit of infinite range.

\paragraph{Paper organization}
Section \ref{sec:basicmodel} introduces all the basic ingredients of the models, including standard RUE and DUE. This is useful for having some reference results for the novel models investigated later on.

Section \ref{sec:V2Vmodel}, which is the core of the paper, introduces two novel models which extend RUE and DUE exploiting V2V communications, and show the results of some numerical tests.

The paper is concluded with some comments and ideas for future research.

\section{The models with full information}\label{sec:basicmodel}
In this section we recall the follow-the-leader microscopic differential model we use for describing the vehicles dynamics. 
Main reference for this part is \cite{cristiani2016NHM}. 
After that, we discuss the process of route choice, following the previously introduced scheme (BB, RUE, DUE).

\paragraph{Notations.}
Let us define a road network $\mathcal N$, composed by $\Nr$ roads indexed by $\r=1,\ldots, \Nr$, and $\Nj$ junctions indexed by $\j=1,\ldots, \Nj$.
Each road $\r$ is unidirectional, and connect a junction $\j^{start}_\r$ (begin) to  another junction $\j^{end}_\r$ (end). Its length is denoted by $L_\r$.
At the same time, each junction $\j$ is associated to a set of incoming roads $\{r^{inc}_\j\}$, and to a set of outgoing roads $\{r^{out}_\j\}$.
This defines completely the geometry of the network.

Let us also introduce notations for the cars moving on $\mathcal N$. 
We consider $\Nc$ cars indexed by $\c=1,\ldots,\Nc$. 
Each car is characterized by the junction $\o_\c$ from where it starts its journey (origin), the junction $\d_\c$ that it wants to reach in minimal time (destination), the road $\r_\c(t)$ on which it is located at any time $t$, the distance $x_\c(t)\in[0,L_{\r_\c}]$ from the origin $\j^{start}_{\r_\c}$ of the current road (coordinate of the segment), and by its speed $v_\c(t)$.

We also denote by $\ell$ the length of any vehicle, and by $\mathcal L$ their total length.
It is important to note that there are two ways for varying the total number of cars $N_\c$: one can do it keeping $\ell$ fixed or one can keep the relation $\ell=\frac{\mathcal L}{\Nc}$ hold true, for a fixed $\mathcal L$.
\REV{In the first setting (actually the one considered in this paper), one truly changes the number of vehicles adding mass to the system, and therefore increases the degree of congestion of the network; in the second setting, instead, one assumes that cars ``shrink'' as their number increases, preserving the total mass of the system. This way, we are just moving from a microscopic description to a macroscopic one, and, in the limit, cars disappear and are substituted by a density function with the same total mass \cite{cristiani2016NHM}.}

Each car $\c$ has a desired path along the network which joins $\o_\c$ and $\d_\c$. 
The path is defined by the driver's choices at junctions, i.e.\ by the function $(t,\j)\to\nextroad_\c(t,\j)$ which denotes the road which should be taken by the car $\c$ if it reaches junction $\j$ at time $t$. 
Clearly $\nextroad_\c(t,\j)\in\{r^{out}_\j\}$.
To enforce the stop at destination, we set $\nextroad_\c(t,\d_\c)=\emptyset$. 
At this point, we assume that the car disappears and cannot be seen by anyone.

Finally, we define the function $t\to\nextcar_\c(t)$ which gives the index of the car preceding car $\c$ at time $t$, \textit{along the path of car $\c$}. 
It is possible that such a car does not exist, in this case car $\c$ is labeled as a \textit{leader} and $\nextcar_\c(t)$ is set to $\emptyset$.

\begin{figure}[ht]
    \centering
    \subfloat[simple 11-road network]{\label{subfig:networks-simple}\includegraphics[width=.48\linewidth]{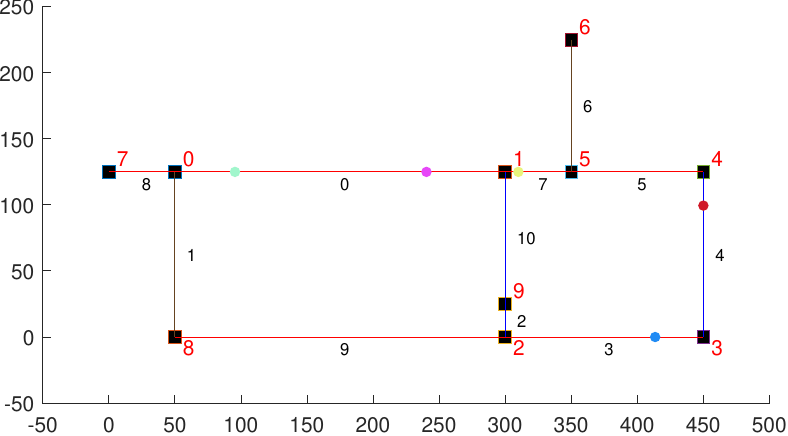}}\quad
    \subfloat[$5\times 5$ Manhattan-like road network with $L_\r=50$]{\label{subfig:networks-manhattan}\includegraphics[width=.48\linewidth]{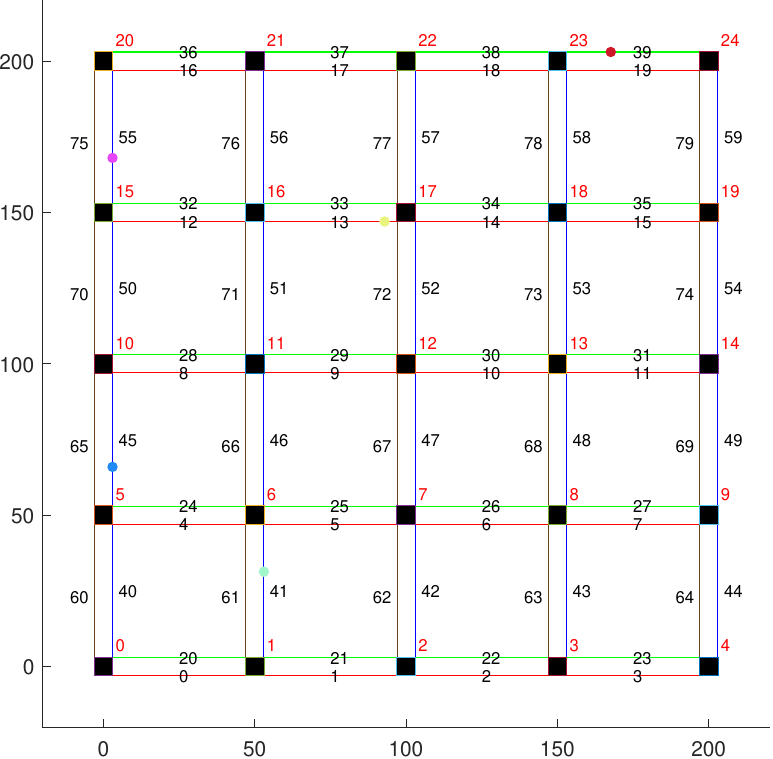}}\par
    \caption{Networks utilized for numerical tests, with road and junction numbering. Roads are all one-way, their colors indicate the direction of motion (blue=up, brown=down, red=right, green=left). Roads between the same pair of junctions are actually superimposed, the shift is only for pictorial purposes. Cars are visualized as colored dots}
    \label{fig:networks}
\end{figure}

\paragraph{Considered networks}
For numerical tests we will consider two road networks, see in Fig.\ \ref{fig:networks}.
The simple 11-road network is mainly used to highlight the difference between RUE and DUE and check the correctness of the numerical code. 
The Manhattan-like network is instead chosen to highlight the role of V2V communications and study the spreading of information in a perfectly symmetric context. 
Note that the second network allows cars to choose many equivalently optimal path to reach their destinations.
\subsection{Dynamic network loading}\label{sec:basicforward}
The vehicles' dynamics is assumed to be first-order, meaning that it can be described by a system of first-order ordinary differential equation
\begin{equation}\label{ODE}
\dot x_\c(t)=v_\c(\delta_\c(t)), \quad t>0, \quad \c=1,\ldots, N_\c,
\end{equation}
where $\delta_\c(t)$ is the distance between car $\c$ and $\nextcar_\c(t)$ at time $t$ ($\delta_\c=+\infty$ for leader cars) and $v_\c(\cdot)$ is a given velocity function. 
The initial condition is given by the couple ($\r_\c(0)=\nextroad_\c(0,\o_\c) ,\ x_\c(0)$). 

We define the velocity function as
\begin{equation}\label{def:vel}
v_\c(\delta):=\left\{\begin{array}{ll}
	0, & \delta<\ell, \\
	\vmax\left(1-\frac{\ell}{\delta}\right), & \delta\geq\ell,
\end{array}
\right. 
\end{equation}
where $\vmax$ is the maximal velocity. \REV{In this way, the velocity is proportional to the headway of the car $\c$, and bounded in the interval $[0,\vmax]$}.
When a car reaches the end of its current road, it continues on the next road along its path as prescribed by the function $\nextroad$. 

\REV{As already noted in \cite{cristiani2016NHM}, the right-hand side in \eqref{ODE} is not continuous in time, since $\delta_\c$ is not; Moreover, this model does not guarantee that the distance between cars is larger than $\ell$ at any time, even if such constraint is satisfied at $t=0$. However, pairs of cars closer than $\ell$ to each other can be found only in small intervals $[L_\r-\ell,L_\r]$ just before any junction, and their number is bounded by the number of incoming roads of that junction, see \cite[Sect.\ 3.1]{cristiani2016NHM} for details. If this happens, the car behind stops and waits to have enough space ahead to move on.}

Finally, as already mentioned in the Introduction, we recall that the micro-to-macro limit of this model was investigated in \cite{cristiani2016NHM}.

\subsection{Static route optimization via Dynamic Programming Principle}\label{sec:basicbackward}
BB just requires to compute the shortest path for any car $\c$ from origin $\o_\c$ to destination $\d_\c$ once and for all, ignoring the presence of other cars on the network. \REV{Since vehicles' dynamics is not involved here, this is a standard discrete optimization problem on a graph.
To solve it}, one can simply employ the Dynamic Programming Principle (DPP).

First of all, we define the static travel time $\TT_s(\r)$ of road $\r$ as
$$
\TT_s(\r)=\frac{L_\r}{\vmax},
$$
which acts as the \textit{cost} (or \textit{weight}) of the road $\r$ in the minimization.
After that, for any car $\c$, we define the value function \REV{$V_\c:\{1,\ldots,N_\j\}\to \mathbb R^+$,} which associates to any junction $\j$ the minimum time to reach the destination $\d_\c$ from that junction along the optimal path.
Note that in this scenario the velocity is constant, then the shortest path coincides with the fastest path.

The DPP states that the value function $V_\c$ is solution to the boundary-value problem
\begin{equation}\label{PPD-BB}
	\left\{
	\begin{array}{l}
	V_\c(\j)=\min\limits_{\r\in\{r^{out}_\j\}} \big\{V_\c(\j^{end}_\r) + \TT_s(\r)\big\},\qquad \forall \j\neq \d_\c, \\ [4mm]
	V_\c(\d_\c)=0,
	\end{array}
\right.  
\end{equation}
\REV{where the set $\{r^{out}_\j\}$, $j\neq \d_\c$, gives the set of admissible controls (in other words, the shortest path is found by choosing, in the optimal way, the outgoing roads from each junction, until $\d_\c$ is reached).}
$V_\c$ can be easily computed by a fixed-point algorithm, i.e.\ iterating the first equation in \eqref{PPD-BB} until convergence, starting from $V_\c(\j)=+\infty$ as initial guess for all $\j\neq \d_\c$.
Note that convergence is always obtained in a \textit{finite} number of steps.

Once $V_\c$ has been computed, the optimal (static) $\nextroad$ function is easily found as 
$$
\nextroad_\c(\j)=\arg\min\limits_{\r\in\{r^{out}_\j\}} \big\{V_\c(\j^{end}_\r) + \TT_s(\r)\big\},\qquad \forall \c, \ \forall\j.
$$

\subsection{Reactive user equilibrium}\label{sec:RUE}
In order to compute the RUE the procedure is similar to the previous case, but it must be repeated at any time $t$ because of the time-dependence of $\TT$.
Indeed, in this case $\TT(t,\r)$ is defined as the time needed to drive along road $\r$, starting from the beginning of the road at time $t$. The travel time depends on the other cars present on the road, thus different congestion levels actually affect the weights $\TT$.

While in the macroscopic setting one can simply define $\TT(t,\r)=\int_\r \frac{dx}{v(\rho(t,x))}$, with $\rho(t,x)$ being the density of vehicles and $v(\rho)$ being their velocity, in a microscopic setting, instead, the computation of $\TT(t,\r)$ is much more tricky.  
Several techniques were used in the literature, and many other papers, especially in the engineering literature, lack some details on this point.
Here are some possible methods:
\begin{enumerate}
	\item[M1.] $\TT(t,\r)$ is the average travel time experienced by the cars passed along the road $\r$ before time $t$, duly weighted in time to prioritize most recent vehicles.
	\item[M2.] $\TT(t,\r)$ is any increasing function of the number of vehicles on the road $\r$ at time $t$.
	\item[M3.] $\TT(t,\r)$ is equal to $L_\r/\hat v_\r$, with $\hat v_\r$ being the average velocity of the cars on the road $\r$ at time $t$.
	\item[M4.] \REV{$\TT(t,\r)$ is the time needed to reach the end of the road $\r$ by a \textit{probe vehicle} that starts moving at time $t$ from the beginning of the road $\r$. The probe vehicle exists only in a fictitious simulation, which is an exact copy of the real one (same network, cars, and dynamics) and it is ran to the side of the main simulation, when needed.}  
\end{enumerate}
In any case, it must be considered that one can get very large values for $\TT$ (or even $+\infty$), because of the presence of queues.
Moreover, yet another difficulty arises when DUE comes into play, we will detail this point in the next section.

The DPP takes the form
\begin{equation}\label{PPD-RUE}
	\left\{
	\begin{array}{l}
		V_\c(\bar t,\j)=\min\limits_{\r\in\{r^{out}_\j\}} \big\{V_\c(\bar t,\j^{end}_\r) + \TT(\bar t,\r)\big\},\qquad \forall \j\neq \d_\c, \\ [4mm]
		V_\c(\bar t,\d_\c)=0,
	\end{array}
	\right.  
\end{equation}
where $\bar t$ is any \textit{fixed} time, which acts as a parameter here in this context.

Then,
$$
\nextroad_\c(\bar t,\j)=\arg\min\limits_{\r\in\{r^{out}_\j\}} \big\{V_\c(\bar t,\j^{end}_\r) + \TT(\bar t,\r)\big\},\qquad \forall \c, \ \forall\j.
$$

\subsection{Dynamic user equilibrium}\label{sec:DUE}
The procedure for computing the DUE is different from the previous ones and more complicated, because the dynamic network loading (\textit{forward-in-time} traffic flow model) and the route optimization (\textit{backward-in-time} DPP) are fully coupled in the space-time and must be computed as one \cite{cacace2019EJAM, camilli2015DCDS, cristiani2015NHM, han2019NSE, morandi20244OR}.

Let us consider a final time $T_{fin}$ large enough to be sure that all cars have reached their destination at that time. The DUE is given by the tuple 
\begin{equation}\label{DUEtuple}
\big(x_\c(t),\ \r_\c(t),\ \TT(t,\r),\ V_\c(t,\j),\ \nextroad_\c(t,\j)\big),
\end{equation}
for all $\c$, $\j$, $\r$, and $t\in[0,\Tfin]$
such that \eqref{ODE}-\eqref{def:vel} are satisfied and
\begin{equation}\label{PPD-DUE}
	\left\{
	\begin{array}{l}
		V_\c(t,\j)=\min\limits_{\r\in\{r^{out}_\j\}} 
		\big\{V_\c(t+\TT(t,\r),\j^{end}_\r) + \TT(t,\r)\big\},\quad \forall \j\neq \d_\c, \quad t\in[0,\Tfin], \\ [4mm]
		V_\c(t,\d_\c)=0, \quad t\in[0,\Tfin], \\ [2mm]
		V_\c(T_{fin},\j)=+\infty, \quad \forall \j\neq \d_\c, \\ [2mm]
		\nextroad_\c(t,\j)=\arg\min\limits_{\r\in\{r^{out}_\j\}} \big\{V_\c(t+\TT(t,\r),\j^{end}_\r) + \TT(t,\r)\big\},\quad \forall\j\neq \d_\c, 
	\end{array}
	\right.  
\end{equation}
is satisfied as well. 
\REV{The function $w$ is again computed by one of the methods M1-M4 mentioned above, and couples the forward and the backward parts of the system of equations.}
Note that it is normal that, at equilibrium, $V_\c=+\infty$ at some $\j$, even for $t<\Tfin$. 
This happens if there is not enough time to reach the destination before $\Tfin$ starting from $\j$ at time $t$. 
The region in the space-time where $V_\c$ is finite is called \textit{reachable set} in control theory \cite[Chapt.\ IV, Sect.\ 1]{bardicapuzzodolcettabook}.

The actual computation of the solution \eqref{DUEtuple} is not easy, and several methods were devised in the literature. 
We have used one of the most common, which consists in iterating the forward-in-time and the backward-in-time equations, keeping fixed $\nextroad_\c$ for all $\c$ in the first case and $\TT$ for all $t$ and $\r$ in the second case.
Iterations keep going until convergence is reached, if any.
To add complexity, many possibilities exist to define convergence in this framework: 
one can check if all positions of all cars at all times stabilize at some iteration, or at least the functions $\nextroad_\c$ stabilize.
We have preferred to employ an even weaker indicator, i.e.\ we check if the \textit{total travel time} $\TTT$, defined as the sum of times needed by all cars to complete their journey,
\begin{equation}\label{def:TTT}
\TTT:=\sum_\c \TTc,\qquad
\TTc:=\min\{t:\j_{\r_\c}^{end}=
\d_\c \text{  and } x_\c=L_{\r_\c}\},
\end{equation}
eventually stabilizes. 
\REV{Doing this, of course, there is the risk of masking different route distributions with the same global $\TTT$.
That said}, in any case the exact convergence is really difficult to obtain in a microscopic framework, and some initial conditions easily lead to instabilities, as reported in many papers; see, e.g., \cite{peeta1995AOR, cheng2013JTRF}. 
To fix this, some techniques were proposed; among them, we have chosen the \textit{method of successive averages} (MSA) (cf.\ also the \textit{fictitious play} technique), consisting in using an average of the values of the weights $w$ obtained in previous iterations rather than just the last available value.

Note that instabilities can also arise at the macroscopic scale, unless a regularization term is added; cf., e.g., \cite{cristiani2015NHM, cristiani2023CMS}.

\begin{remark}[Match between forward and backward dynamics]
It is useful to note here that the link between the forward-in-time model \eqref{ODE}-\eqref{def:vel} and the backward-in-time model \eqref{PPD-DUE} is obtained by means of the function $\TT$. In fact, drivers optimize their path assuming a certain travel time for each road at each time, and then that travel time corresponds exactly to that one experienced by a car which passes on that road at that time. 
In principle, there is no reason to require that the two travel times (the real one and the one used as weight in the optimization process) coincide, because one can assume that drivers do not know exactly the real traffic flow dynamics because of a wrong perception.
At the same time, one can expect that having perfect symmetry in forward and backward dynamics can reduce instabilities of DUE.

\REV{Method M1 is robust and temporally smooth but reacts slowly to rapid congestion changes, since it relies on past observations, therefore it may delay path updates. Moreover, it needs at least an additional parameter for prioritizing most recent vehicles.}

\REV{Method M2 is the simplest and most economical one. However, it only reflects instantaneous density and does not distinguish different headway configurations, queues or stop-and-go patterns, making it rather coarse.}

\REV{Method M3 is easy to compute, reacts quickly to current congestion, and is almost aligned with the microscopic dynamics, since speeds directly encode headway.}

\REV{Finally, method M4 seems to guarantee perfect symmetry in forward and backward dynamics; also, it is in line with the DPP in the space-time \eqref{PPD-DUE}, which states that if one moves along the optimal path (unknown), ``the time needed to reach the destination from any junction $\j$ at time $t$ equals the time $\delta t$ to reach the following junction $\j^*$ plus the time to reach the destination starting from junction $\j^*$ at time $t+\delta t$''. The time $\delta t$ is indeed the time computed in the fictitious simulation. Actually, we observed that method M4 reduces instabilities for some initial conditions, but the improvement is not always consistent.} 
	
\REV{Considering the major difficulty in implementing M4, and the increased computational time, we think that the best balanced choice among all is M3; therefore, we adopt it in all our simulations.}  
\end{remark}


%
%
%
%
%
%

\subsection{Numerical tests}
To conclude this section, we present some preliminary numerical tests which will serve as a reference for the results of the novel model detailed in Section \ref{sec:V2Vmodel}.
Equation \eqref{ODE} is solved by means of the Explicit Euler scheme with time step $\Delta t$. 
The model parameters are chosen as summarized in Table \ref{tab:parameters}.
For computing the weights $w$ we have always used the method M3, but DUE was also computed using method M4 for comparison.
The final time $\Tfin$ was always chosen large enough so that all cars eventually reach their destination.
The cars' length $\ell$ in \eqref{def:vel} includes the minimal free space usually left between cars, even if full congestion is reached.
\begin{table}[h!]
	\centering
	\begin{tabular}{l|l|l|l}
	Parameter & Symbol & Value &Unit \\ \hline
	Time step                     &$\Dt$        &0.6             &sec   \\
	Number of junctions (Manhattan-like netw.)  &$N_\j$       &$3^2,5^2,7^2$   & --   \\
	Length of each road (Manhattan-like netw.)  &$L_\r$       &50 -- 300       & m    \\   
	Number of vehicles            &$N_\c$       &25 -- 100       & --   \\
	Maximal velocity              &$v_{max}$    &50              & km/h \\
	Length of each vehicle        &$\ell$       &10              &m     \\
	\end{tabular}
	\caption{Parameters used in numerical simulations}
	\label{tab:parameters}
\end{table}

\subsubsection{Test 0: Basic statistics} 
\begin{figure}[h!]
	\centering
	\subfloat[Comparison for $N_\c=25$, $50$, and $100$ cars on $5 \times 5$ Manhattan-like network]{\label{subfig:mediacumulata_Nc}
		\includegraphics[width=.47\linewidth]{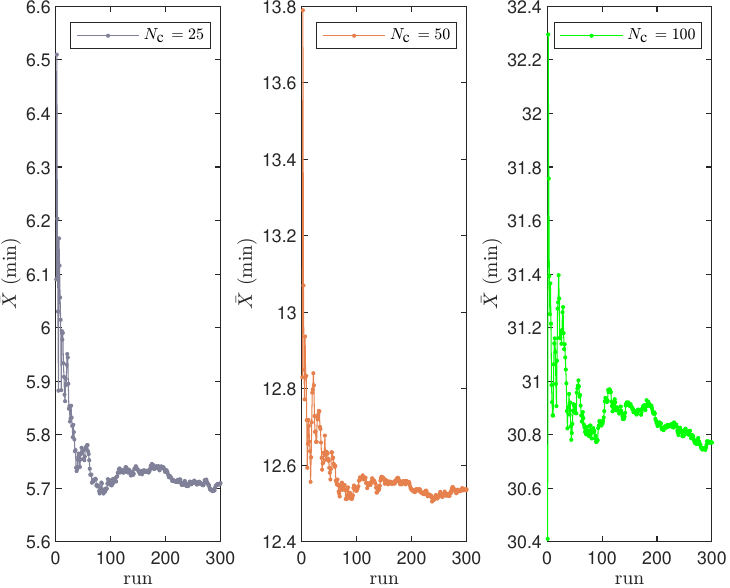}}\quad
	\subfloat[Comparison for $3\times 3$, $5\times 5$, and $7\times 7$ Manhattan-like network with $N_\c=100$ cars]{\label{subfig:mediacumulata_Nj}
		\includegraphics[width=.47\linewidth]{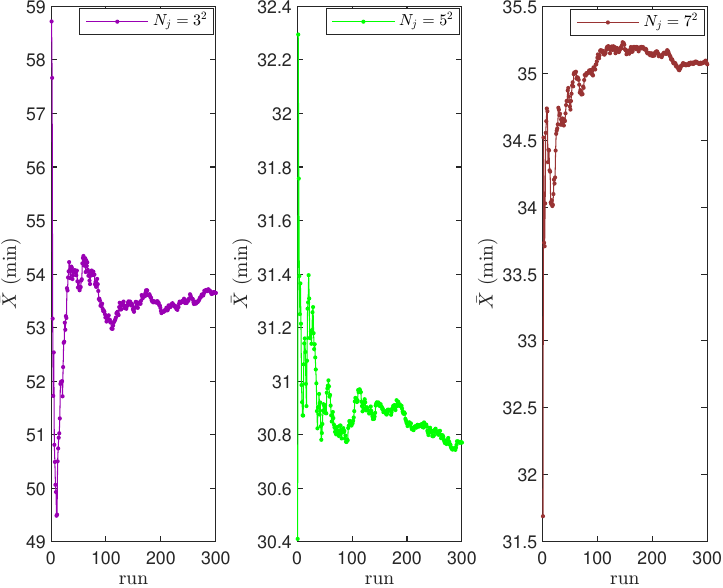}}\par
	\caption{Test 0. Cumulative average of $\TTT$ obtained for RUE, with $L_\r=50$ and random OD pairs}
	\label{fig:mediacumulata}
\end{figure}

Many tests we will perform in the following involve random origins and/or destinations.
This is also done because we want to keep an easy-to-get macro\-scale counterpart of our experiments, and random uniform initial and final conditions are associated to a constant initial and final densities on all over the network.
Working with random variables requires to observe average quantities (expected values), therefore the question arises how many runs (repetitions) have to be performed in order to have meaningful statistics. 
The quantity which we will study the most is the total travel time $\TTT$ defined in \eqref{def:TTT}.
\REV{Moreover, we expect RUE be the dynamics with larger deviations, so we run 300 simulations for that behavior on the Manhattan-like network with random OD pairs, $L_\r=50$, and we plot the cumulative average}
$$
\bar X_{r}:=\frac{1}{r}\sum_{i=1}^{r}\TTT(i), \qquad r=1,\ldots,300,
$$
where $\TTT(i)$ is the total travel time for run $i$; see Fig.\ \ref{fig:mediacumulata}.
After that, we use the standard formula
\begin{equation}
	\mathbb P \left( \text{\textbar} \bar{X}_r-\mu \text{\textbar} \leq \varepsilon_{r,\alpha}\right)=1-\alpha, \qquad
	\varepsilon_{r,\alpha}:=t_{\frac{\alpha}{2},r-1}\frac{\bar\sigma_r}{\sqrt{r}},
\end{equation}
where $\mu$ is the true mean of $\TTT$ (approximated by $\TTT(500)$), $\bar\sigma_r$ is the standard deviation of the $r$ values, $t_{\frac{\alpha}{2},r-1}$ is the Student's $t$ distribution with $r-1$ degrees of freedom, and $\alpha\in[0,1]$ is a parameter related to the reliability of the estimation.
By setting $\varepsilon_{r,\alpha}=0.25$ minutes and $\alpha=0.01$, we find that performing 300 runs is a reasonable choice for RUE, then it is for BB and DUE too. 
Consequently, in the tests that follow we will always provide results obtained by averaging over 300 runs.

\subsubsection{Test 1: BB vs. RUE vs. DUE} 
\begin{figure}[t!]
    \centering
    \subfloat[BB (\href{www.emilianocristiani.it/attach/paper_trafficV2Vmicro_smallnetwork_50cars_BB.mp4}{video})]{\label{subfig:sim_simple_netw-BB}
    \includegraphics[width=.45\linewidth]{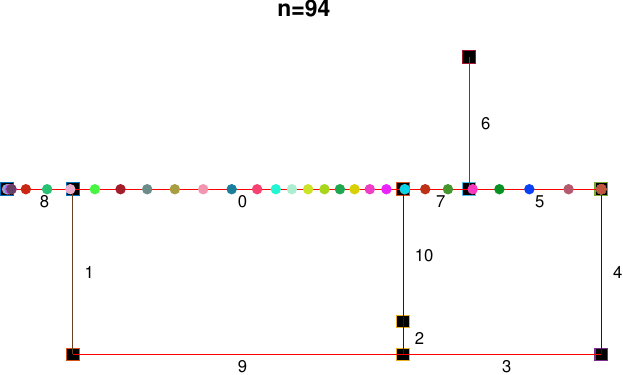}}\\
    \subfloat[RUE (\href{www.emilianocristiani.it/attach/paper_trafficV2Vmicro_smallnetwork_50cars_RUE.mp4}{video})]{\label{subfig:sim_simple_netw-RUE}
    \includegraphics[width=.45\linewidth]{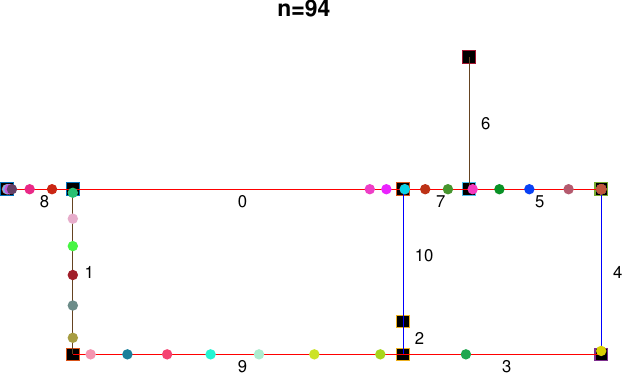}}\quad
    \subfloat[DUE (\href{www.emilianocristiani.it/attach/paper_trafficV2Vmicro_smallnetwork_50cars_DUE.mp4}{video})]{\label{subfig:sim_simple_netw-DUE}
    \includegraphics[width=.45\linewidth]{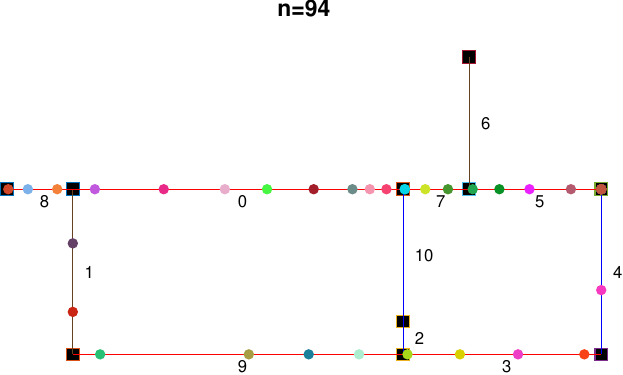}}
    \caption{Test 1. Three screenshots, corresponding to the same time $t$, of the simulation on the 11-road simple network (colors are just pictorial). They clearly show the difference among the three degrees of rationality}
    \label{fig:sim_simple_netw}
\end{figure}

In this test we compare BB, RUE, and DUE on the two networks depicted in Fig.\ \ref{fig:networks}.
Let us start with the simple 11-road network, which is specifically designed to highlight the differences among the three degrees of rationality.
The destination is junction 4 for all $N_\c=50$ cars, while initial positions are distributed across roads 2, 6, and 8.
BB prescribes that all cars follow the shortest path, in particular cars starting from road 8 move straight through roads 0, 7, and finally 5 (Fig.\ \ref{subfig:sim_simple_netw-BB}).\footnote{\scriptsize\REV{\detokenize{www.emilianocristiani.it/attach/paper_trafficV2Vmicro_smallnetwork_50cars_BB.mp4}}}
RUE instead prescribes that cars starting from road 8 initially choose the shortest path, then they switch to the longer path 1-9-3-4 whenever the shortest becomes congested due to the merge with the cars coming from 2 and 6. When also the longer path becomes congested as well, they re-switch to the shortest one (Fig.\ \ref{subfig:sim_simple_netw-RUE}).\footnote{\scriptsize\REV{\detokenize{www.emilianocristiani.it/attach/paper_trafficV2Vmicro_smallnetwork_50cars_RUE.mp4}}}
Finally, DUE prescribes that cars starting form 8 are able to forecast the congestion on road 7 well before reaching the first junction, therefore they immediately start alternating between roads 0 and 1 in such a way that the shortest and the longer paths become competitive with each other (Fig.\ \ref{subfig:sim_simple_netw-DUE}).\footnote{\scriptsize\REV{\detokenize{www.emilianocristiani.it/attach/paper_trafficV2Vmicro_smallnetwork_50cars_DUE.mp4}}} 

\begin{figure}[b!]
    \centering
    \subfloat[11-road network with fixed OD pairs]{\label{subfig:BB_vs_RUE_vs_DUE_fixed}
    \includegraphics[width=.45\linewidth]{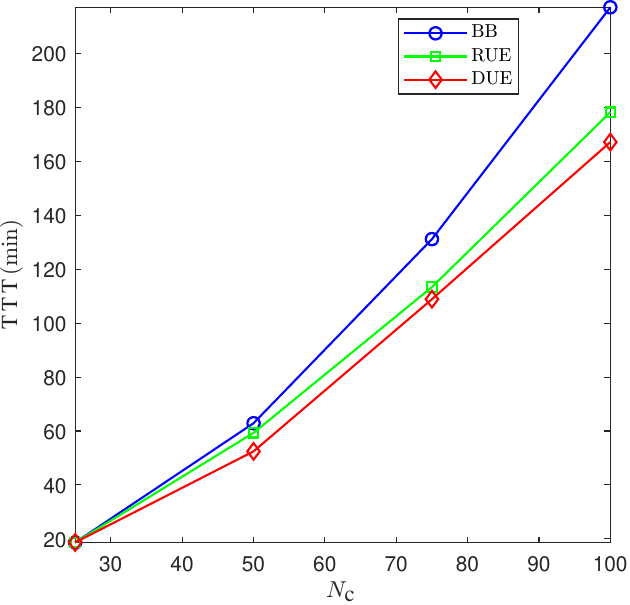}}\quad
    \subfloat[Manhattan-like $5\times 5$ network with random OD pairs and $L_\r=50$]{\label{subfig:BB_vs_RUE_vs_DUE_random}
    \includegraphics[width=.45\linewidth]{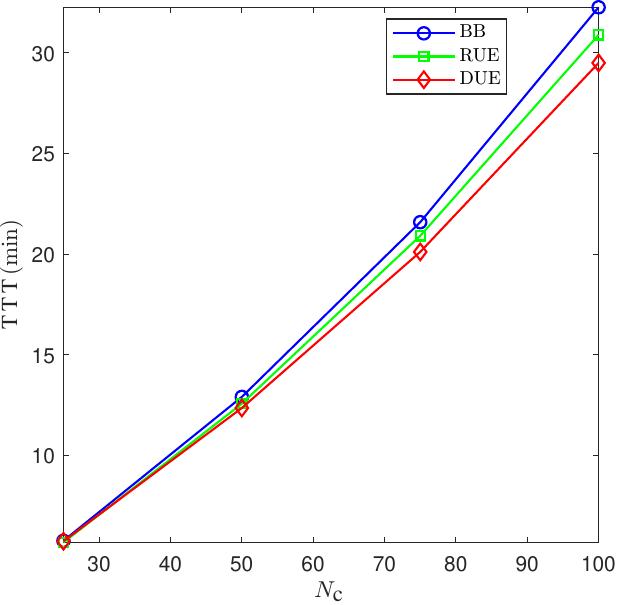}}\par
    \caption{Test 1. Comparison between average $\TTT$'s obtained with 25, 50, 75, and 100 cars and different degrees of rationality (BB, RUE, DUE)}
    \label{fig:BB_vs_RUE_vs_DUE}
\end{figure}

The difference among the three degrees of rationality is clearly visible also looking at the $\TTT$, see Fig.\ \ref{subfig:BB_vs_RUE_vs_DUE_fixed}.
As expected, DUE performs better than RUE which, in turn, performs better than BB. 
In addition, the improvement due to the higher degree of rationality increases with the number of cars.

%
%
%

We have also run a similar test on the $5\times 5$ Manhattan-like network with random OD pairs and $L_\r=50$. 
In this case, the three behaviors are more similar to each other in comparison to the previous case. 
This can be well observed by comparing the $\TTT$s for different values of $N_\c$, see Fig.\ \ref{subfig:BB_vs_RUE_vs_DUE_random}.
We think that the main reasons for that are the large number of equivalently optimal path and the less degree of congestion.

%
%
%
%
%
%

\section{The models with V2V communications}\label{sec:V2Vmodel}
In this section we introduce the main novelty of the paper, namely V2V communications. 
The idea is that any time two vehicles come closer than a threshold distance $\radius$ to each other, they are able to share information.
In particular, we assume that any car $\c$ can possibly share:
\begin{itemize}
	\item its current position (road $\r_\c$ and coordinate $x_\c$); 
	\item its destination $\d_\c$;
	\item the whole planned path toward destination, i.e.\ the function $\nextroad_\c$.
\end{itemize}
All of these pieces of information come along the \textit{timestamp} of the rendez-vous and are stored in the OPUs of both vehicles.
We also assume that vehicles communicate every $\Dcom$ time units (so as not to burden the radio network) and that vehicles \textit{forget} the information older than $\Dmem$ time units. This introduces a \textit{memory} in the system.
Lastly, the \textit{cascade effect}: if enabled, vehicles transmit not only the information related to their own, but also \textit{all the pieces of information related to the other vehicles} they have collected and that they have not forgotten yet. 

These additional parameters dedicated to V2V models are summarized in Table \ref{tab:additionalparameters}.
\begin{table}[h!]
	\centering
	\begin{tabular}{l|l|l|l|l}
		Parameter            & Symbol     & Reference value & Values explored  &Unit \\ \hline
		Communication range  &$\radius$   & 150              & $[0,+\infty)$    & m        \\
            Communication pause	 &$\Dcom$     & 0                & $[0,30]$         & sec      \\ 
		Memory               &$\Dmem$     & $+\infty$        & $[0, +\infty)$   & min      \\
		Cascade effect       & --         & yes              &yes/no            & --       
	\end{tabular}
	\caption{Additional model parameters dedicated to V2V communications}
	\label{tab:additionalparameters}
\end{table}

\subsection{Test 2: spreading of information}\label{sec:spreadingV2V}
Before investigating the role of the various degrees of rationality in the context of V2V communications, we want to focus on the spreading of information across the network. 

We consider the $5\times 5$ Manhattan-like network with random OD pairs and $L_\r=300$. 
The cars start to move at the same time, and we label them as ``active'' until they reach their destination. 
At that point, they are labeled as ``inactive''.
In addition, we denote by $N_a(t)$ the number of active cars at time $t$ (always bounded by $N_\c$), by $\mathcal K_\c(t)$ the set of active cars known by car $\c$ at time $t$ (it is well defined only if $\c$ is active), and by $K_\c(t)$ its cardinality.
The quantity of our interest will be then
\begin{equation}\label{def:KNOW}
\KNOW(t):=\frac{1}{N_a(t)}\sum_{\text{$\c$ is active}} K_\c(t),
\end{equation}
which expresses a sort of global knowledge across the network. 
If $\KNOW(t)=N_a(t)$, the knowledge saturates, since all active cars know all other active cars.

\begin{figure}[t!]
    \centering
    \subfloat[Cascade yes/no]{\label{subfig:spreading-cascade}
    \includegraphics[width=.41\linewidth]{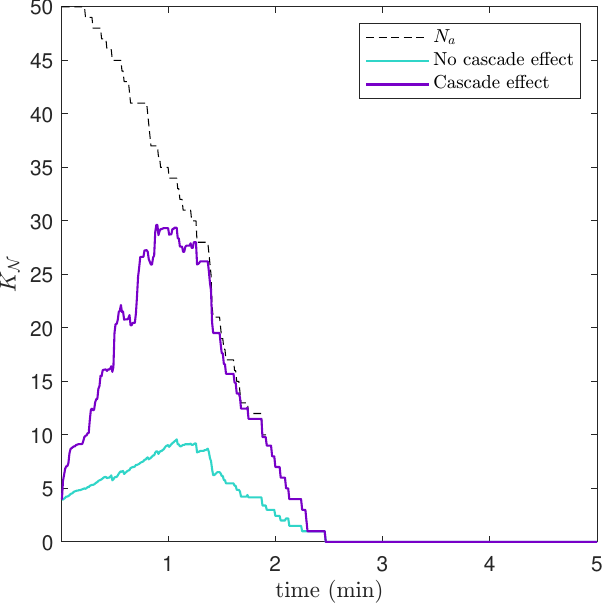}}\quad
    \subfloat[Varying $\Dcom$]{\label{subfig:spreading-dc}
    \includegraphics[width=.41\linewidth]{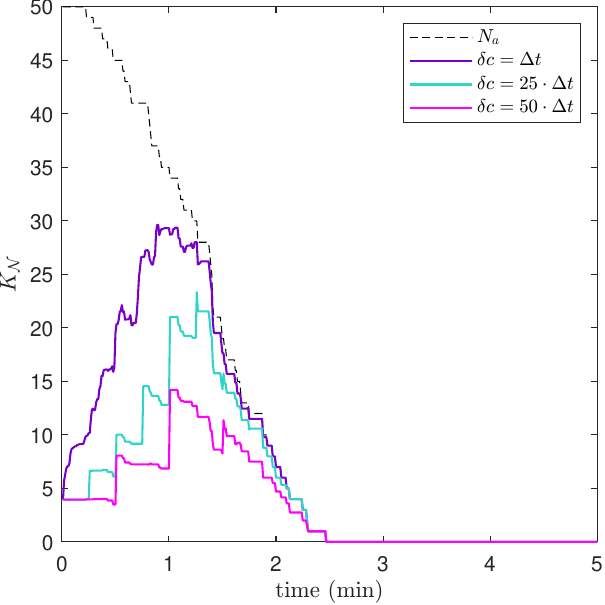}}\\
    \subfloat[Varying $\radius$]{\label{subfig:spreading-R}
    \includegraphics[width=.41\linewidth]{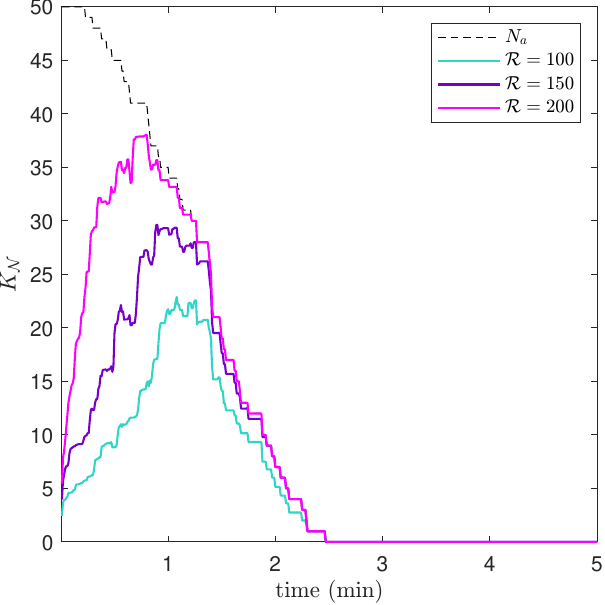}}\quad
    \subfloat[Varying $N_\c$]{\label{subfig:spreading-Nc}
    \includegraphics[width=.41\linewidth]{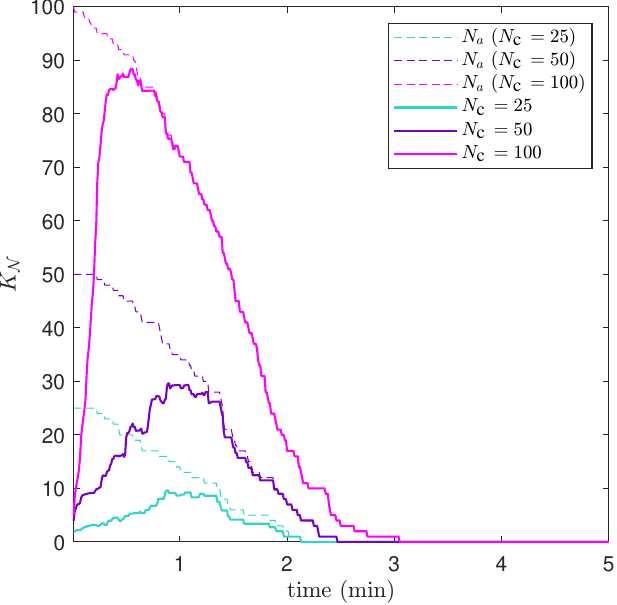}}\\
    \subfloat[Varying $N_\j$]{\label{subfig:spreading-Nj}
    \includegraphics[width=.41\linewidth]{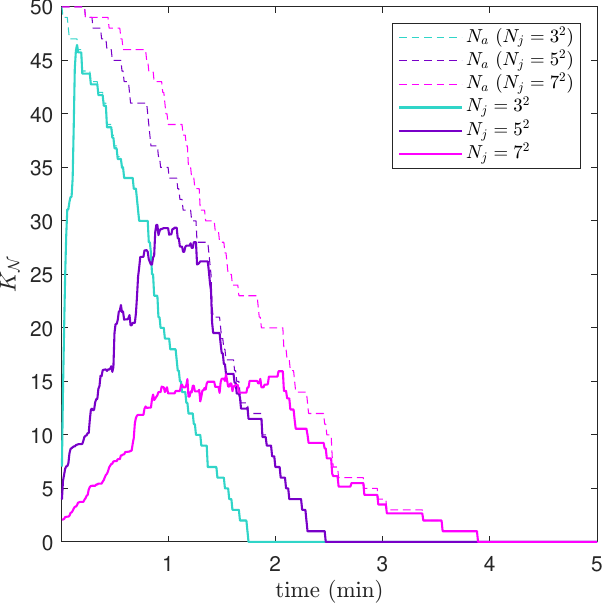}}\quad
    \subfloat[Varying $\Dmem$]{\label{subfig:spreading-dm}
    \includegraphics[width=.41\linewidth]{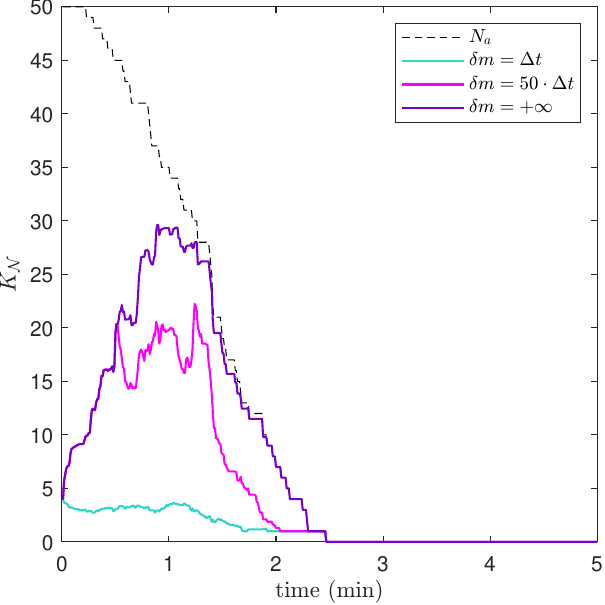}}\par
    \caption{Test 2. Spreading of information}
    \label{fig:spreading}
\end{figure}

In order to simplify the numerical tests, in the following we limit ourselves to the BB, i.e.\ we assume that all cars follow the (static) shortest path to their destination.
Moreover, we use reference parameters as given in Table \ref{tab:additionalparameters}, 
and then we increase/decrease them one at a time in order to study the sensitivity of $\KNOW$ to those parameters; see Fig.\ \ref{fig:spreading}. 
For an easier comparison, in all figures the result obtained with the reference parameters is repeated (purple line).
In addition, we also plot $N_a(t)$ as an upper bound for the knowledge indicator $\KNOW(t)$.

Fig.\ \ref{subfig:spreading-cascade} shows the impact of the cascade effect. As expected, removing the cascade effect slows down the spreading of the information. Also, disabling the cascade effect, not only $\KNOW$ is lower, but \textit{it never saturates}.
Fig.\ \ref{subfig:spreading-dc} shows the sensitivity with respect to the parameter $\Dcom$.
Peaks in correspondence of communications are clearly visible. 
In between, $\KNOW$ remains constant or decreases as some cars become inactive.
Fig.\ \ref{subfig:spreading-R} shows the sensitivity with respect to the parameter $\radius$. 
One can see that if $\radius$ is small, $\KNOW$ does not saturate even if the cascade effect is enabled.
Fig.\ \ref{subfig:spreading-Nc} shows the sensitivity with respect to the parameter $N_\c$.
As expected, $\KNOW$ is proportional to $N_\c$, and the first time $\KNOW$ saturates is inversely proportional to it.
This parameter seems to be the one with the most impact on the information spreading.
Fig.\ \ref{subfig:spreading-Nj} shows the sensitivity with respect to the parameter $N_\j$.
In this case we observe the opposite behavior with respect to the previous case, i.e.\ $\KNOW$ is inversely proportional to $N_\j$, while first time of saturation is proportional to it.
Finally, Fig.\ \ref{subfig:spreading-dm} shows the sensitivity with respect to the parameter $\delta m$.
If cars have no memory the number of known cars equals the cars within the distance $\radius$ at any given time.
For small or mid-size values of $\Dmem$, the function $\KNOW$ does not saturate.

\subsection{Reactive user equilibrium with V2V communications}\label{sec:RUE-V2V}
In this section, we propose a model for investigating the RUE in the context of V2V communications (V2V-RUE).
The idea is that any rational decision-making (path choice, in our case) must be limited by the knowledge acquired by vehicles through communications.

As recalled in Section \ref{sec:RUE}, the RUE is based on the assumption that vehicles cannot make any prediction of the others' future positions but continuously re-optimize their path on the basis of the knowledge of the current state of the whole system.
The main difficulty we face here is the fact that rendez-vous happen at \textit{different times}, meaning that each car has a knowledge of others' positions which is sparse in time. 
This means that even if car $\c$ has encountered many cars in the past, it has no knowledge of the \textit{current} state of the system. 

Before giving the details of the new model, we need to modify the weight function $w$, adding a dependence on the car $\c$ itself, i.e.\ let us consider the function $(\c,\bar t,\r)\to w_\c(\bar t,\r)$. 
This is a crucial step because the weight of a road is computed observing cars staying or moving on that road, and, in our context, each car $\c$ has a different knowledge of the traffic.
More precisely, each car is aware only of the subset $\mathcal K_\c$ of cars actually on the roads, therefore each car can do, in general, a different estimate of the road weights.

\subsubsection{The V2V-RUE model}
The model we propose is based on two main assumptions: 
1) at each rendez-vous, cars share their current position and destination only; 
2) at each fixed time $\bar t$, each car is able to \textit{nowcast}\footnote{
    \textit{nowcast} is a term often used in the engineering literature which merges the words `now' and `forecast', meaning a prediction from the past to the current time (now).
} 
the positions of the known cars by means of a parallel fictitious simulation.
In the parallel simulation done by car $\c$, we assume that the network is populated only by cars in $\mathcal K_\c$, and that 
they move along a minimum-time path, computed, in turn, employing the road weights $w_\c(\bar t,\cdot)$, i.e.\ the weights based on the $\c$'s knowledge only. 

Summarizing, at any fixed time $\bar t$ and for any car $\c$:
\begin{enumerate}
	\item car $\c$ forgets all known cars encountered more than $\delta m$ time units in the past;
	\item positions of all cars $\c^\prime\in\mathcal K_\c(\bar t)$ are nowcast, \REV{in an imaginary} fictitious simulation, from the last rendez-vous time $\hat t$ to time $\bar t$ using \eqref{ODE}-\eqref{def:vel}, assuming that only cars in $\mathcal K_\c(\bar t)$ are present on the network. For doing this, a temporary $\nextroad^\sharp_{\c^\prime}(s,\cdot)$, $s\in[\hat t,\bar t)$ is computed by \eqref{PPD-RUE} using $w_\c(s,\cdot)$, $s\in[\hat t,\bar t)$; 
	\item if $\delta c$ time units have passed since the last communications, car $\c$ checks for new rendez-vous in order to update $\mathcal K_\c(\bar t)$. 
	It gets current position and destination of any other car $\c^\prime$ within communication distance $\radius$. 
	If car $\c^\prime$ was already previously encountered, new pieces of information overwrite the old ones;
	\item weights $w_\c(\bar t,\cdot)$ are computed using cars in $\mathcal K_\c(\bar t)$ only; 
	\item $\nextroad_\c(\bar t,\cdot)$ is computed using \eqref{PPD-RUE};
        \item $\nextroad_\c(\bar t,\cdot)$ is used in \eqref{ODE}-\eqref{def:vel} to move car $\c$ ahead.
\end{enumerate}

\begin{remark}[Nowcast inaccuracy]\label{rem:discrepancy}
    The nowcast procedure does not provide, in general, the correct estimation of the current position of the known cars. 
    This happens because in Step 2, a car $\c^\prime$ known by $\c$ is imaginarily moved along the network using $\nextroad^\sharp_{\c^\prime}$ rather than $\nextroad_{\c^\prime}$. 
    This generates a discrepancy between the estimate and the real positions of known cars of any car, i.e.\ between the factors that influence the optimal path planning and those that influence the real vehicle dynamics. 
    The discrepancy disappears only for large enough values of $\radius$, since in that case all cars get full knowledge of the system. 
\end{remark}

\subsubsection{Test 3: the ``blissful ignorance'' paradigm}\label{sec:blissfulignorance}
In this test we consider both $3\times 3$ and $5\times 5$ Manhattan-like road network, with $L_\r=50$, $N_\c=100$, $\delta c=0$ (actually corresponding to $\delta c=\Delta t$ in the numerical discretization, i.e.\ rendez-vous are continuously enabled), $\Dmem=+\infty$ and no cascade effect. 
Regarding OD pairs, we consider:
\begin{itemize}
    \item random origins and destinations;
    \item random origins and a single fixed destination (at the top-center of the network);
    \item fixed origins (spread on two roads) and a single fixed destination (at the top-right corner).
\end{itemize}

In this test we study in particular the behavior of the $\TTT$ for the V2V-RUE as a function of $\radius$. 
The two extreme cases are already known: for $\radius=0$ we recover the BB (since all cars behave as they were alone on the network), while for $\radius\to +\infty$ we recover the standard RUE (since all cars are aware of the exact and up-to-date position of any other car at any current time).

For intermediate values of $\radius$ a nice feature comes to light. 
In fact, one could expect that the larger the communication range, the more information is available, and the lower is $\TTT$. 
Surprisingly enough, this is not always true, since the lowest values of $\TTT$ can be obtained, in some cases, with mid-size communication ranges.
This happens, e.g., in the case of random destinations (Fig.\ \ref{subfig:TTT(R)-V2V-RUE-DUE-3x3-random}-\ref{subfig:TTT(R)-V2V-RUE-DUE-5x5-random}), but not if we consider a unique destination for all cars (at top-right corner, see Fig.\ \ref{subfig:TTT(R)-V2V-RUE-DUE-3x3-fixed}, or at the top-center of the network, see Fig.\ \ref{subfig:TTT(R)-V2V-RUE-DUE-5x5-fixedrandom}).

In our opinion, this ``blissful ignorance'' paradigm can be explained by considering that large values of $\radius$ require each car $\c$ to take into account many other cars which are far from $\c$, and probably will never meet $\c$. If $\c$ deviates from shortest path to avoid such cars, it will probably do something useless, if not disadvantageous. In contrast, if $\c$ knows only nearby cars, it will take into account cars which are more likely to be really encountered along the trip. 
Also, it seems that ignoring far vehicles can be so advantageous in some cases that it compensates the fact that a small communication range tends to increase the error in estimating the positions of the known cars (Remark \ref{rem:discrepancy}).

\subsubsection{Test 4: a new traffic equilibrium}\label{sec:equikibrioignoto}
An intriguing aspect of the V2V-RUE is the possibility of defining a new type of \textit{system equilibrium}. 
Let us explain why we expect it exists by means of a specific example: consider a network populated by many cars, and focus on a car $\c$ in particular with a given OD pair and a mid-size communication range $\radius$. 
Assume that the shortest path $P_s$ joining origin and destination is highly congested, so that it does not coincide with the fastest path.
Assume also that a longer -- and less congested -- path exists, say $P_f$, and that it is actually the fastest path among all possible ones (the one to be preferred).

The interesting point is that neither path can actually be chosen by car $\c$. 
If we force the car to choose the shortest $\&$ congested path $P_s$, the car will encounter a large number of other cars (since the path is congested), then the information collected will allow the car to understand that a longer but less congested path is preferable.
Conversely, if we force the car to choose the longer $\&$ uncongested path $P_f$, the limited amount of traffic-related information will not allow the car to recognize that path as optimal, instead the shortest path will be preferred (if car $\c$ is not aware of any congestion, it will naturally tend to prefer the shortest path). 

In conclusion, we expect there exists a third path which represents a balance point between congestion degree, amount of traffic information, and travel time. 
If it exists, \textit{the equilibrium path $P^*_\c$ for car $\c$ will be characterized by the fact that it allows car $\c$ to encounter the right amount of cars, and then to have the right amount of traffic-related information, such that the path $P^*_\c$ itself is recognized to be the fastest path among all possible ones.}
 
\begin{figure}[ht!]
    \centering
    \subfloat[initial condition]{
    \label{subfig:equilibrium_ic}
    \includegraphics[width=.4\linewidth]{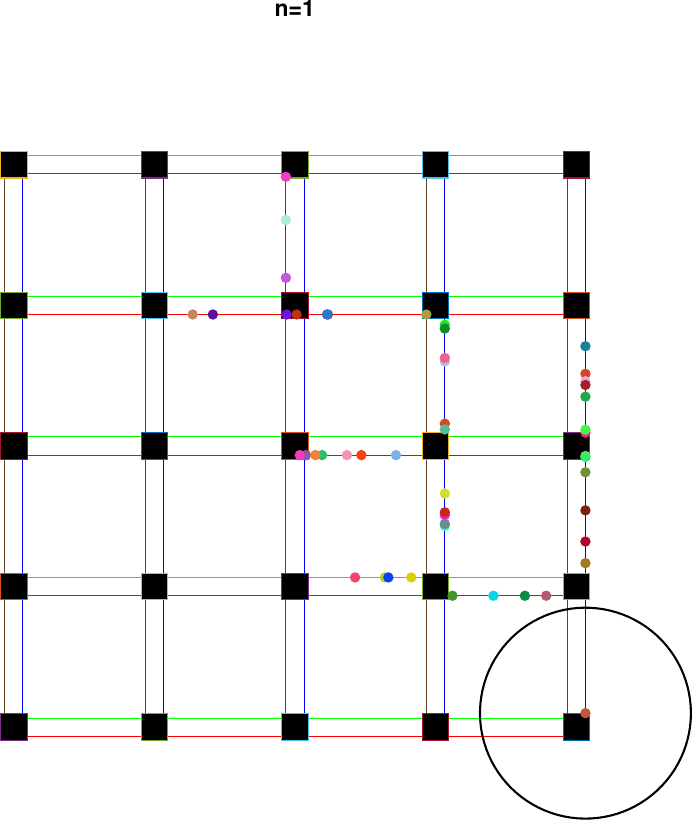}}\\
    \subfloat[shortest path (\href{www.emilianocristiani.it/attach/paper_trafficV2Vmicro_shortest.mp4}{video})]{
    \label{subfig:equilibrium_sp}
    \includegraphics[width=.30\linewidth]{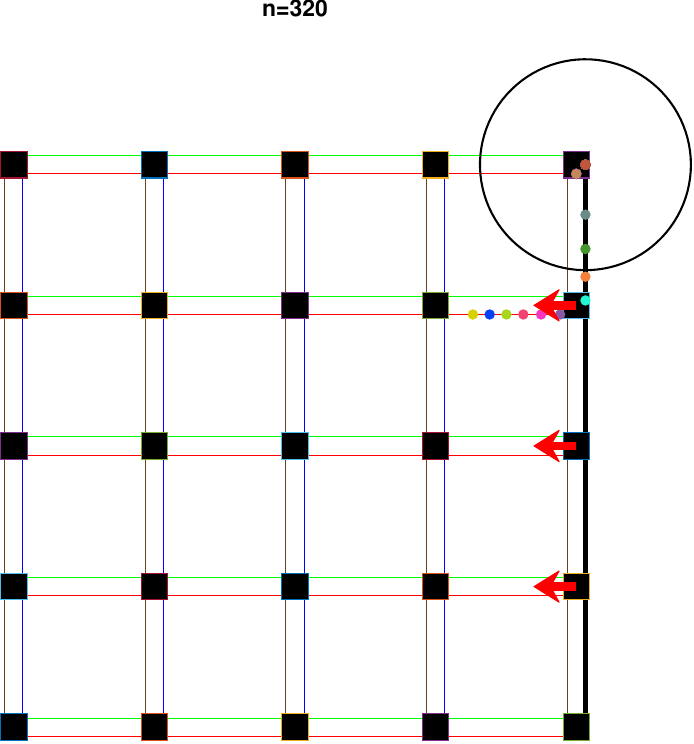}}\quad 
    \subfloat[optimal path (\href{www.emilianocristiani.it/attach/paper_trafficV2Vmicro_fastest.mp4}{video})]{
    \label{subfig:equilibrium_op}
    \includegraphics[width=.30\linewidth]{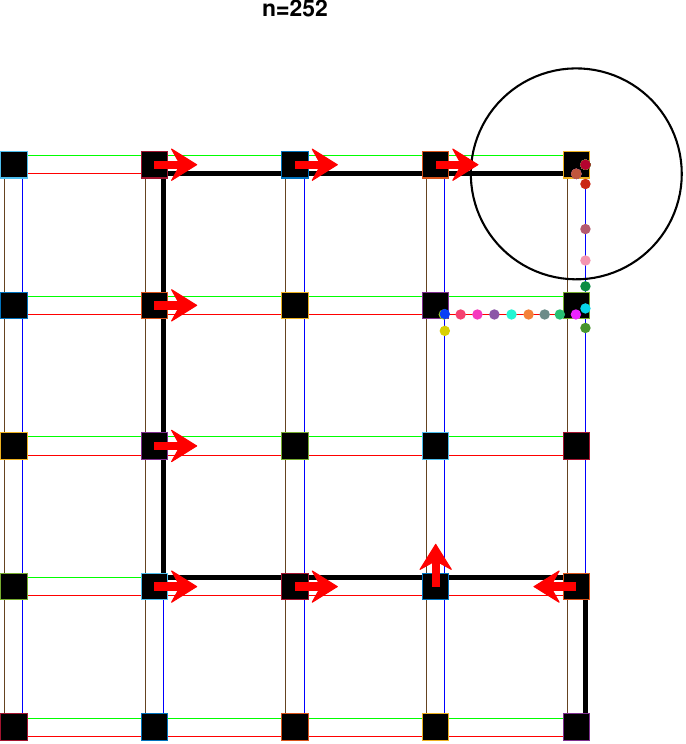}}\quad
    \subfloat[equilibrium path (\href{www.emilianocristiani.it/attach/paper_trafficV2Vmicro_equilibrium.mp4}{video})]{
    \label{subfig:equilibrium_ep}
    \includegraphics[width=.30\linewidth]{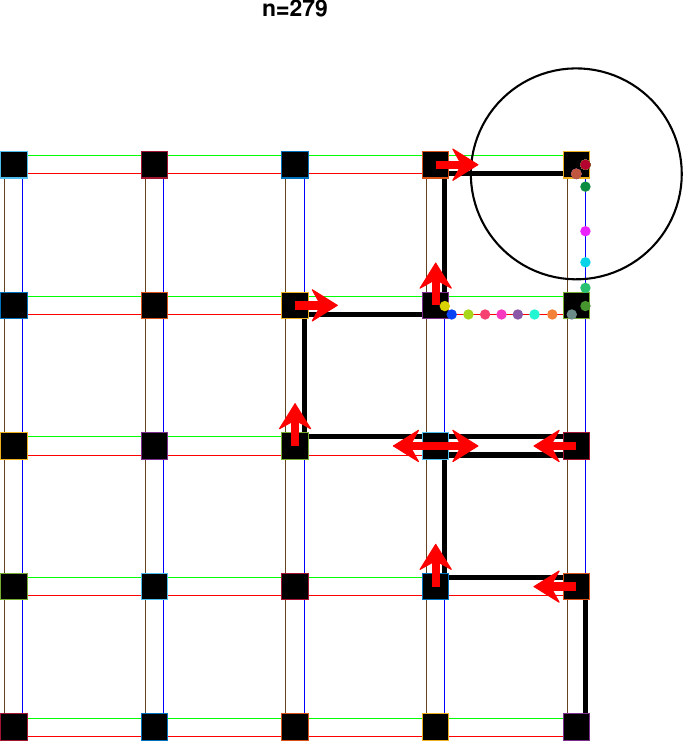}}\par
    \caption{Test 4: A system configuration such that an equilibrium path exists. (a) The tracked car starts from the bottom-right corner and must reach the top-right corner. 
    (b) The shortest path (thick black line) corresponds to a straight line in vertical direction, but it is very congested. Arrows show that V2V-RUE dynamics would suggest to move leftward and choose a longer path. 
    (c) The optimal path corresponds to a much longer path with respect to the shortest one (thick black line). This path is free of congestion. Arrows would suggest to move towards a shorter path.
    (d) The equilibrium path (thick black line) is quite contorted and curvy, but it perfectly matches the red arrows}
    \label{fig:equilibrium}
\end{figure}

Fig.\ \ref{fig:equilibrium}\footnote{
	\scriptsize\REV{Fig.\ \ref{fig:equilibrium}b: \detokenize{www.emilianocristiani.it/attach/paper_trafficV2Vmicro_shortest.mp4}}, \\
	\scriptsize\REV{Fig.\ \ref{fig:equilibrium}c: \detokenize{www.emilianocristiani.it/attach/paper_trafficV2Vmicro_fastest.mp4}}, \\
	\scriptsize\REV{Fig.\ \ref{fig:equilibrium}d: \detokenize{www.emilianocristiani.it/attach/paper_trafficV2Vmicro_equilibrium.mp4}}
} 
shows a specific system configuration on the $5\times5$ Manhattan-like network with $N_\c=50$, for which these three paths (shortest, longer optimal, equilibrium) can be explicitly computed. 
The simulation is realized by manually assigning origin, destination and whole path to each car on the network.
In addition, only for the car of interest, 
we depict at each passed junction a red arrow indicating the choice that the driver would have made if he had been free to decide his path with V2V-RUE dynamics.
If all arrows point toward the road actually chosen, the path is of equilibrium.

\subsection{Dynamic user equilibrium with V2V communications}\label{sec:DUE-V2V}
A model for computing the DUE using V2V communications can be devised in several ways, but, in our opinion, the necessary ingredient should be the possibility, for any car, to \textit{forecast} the path of some other cars, in order to optimize its own path on the basis of that prediction. 
Also, in a V2V context, it is natural to assume that prediction can be made only for known (encountered) cars.

Along these lines, a natural generalization can be derived by assuming that each car $\c$, at any given time $t$, finds its local Wardrop equilibrium with the cars it knows. 
This means that car $\c$ plays a solo game which consists in iterating forward-in-time and backward-in-time equations involving only itself and cars it knows. 
Note that this should be all done with the information possessed by $\c$, and should be repeated at any time because the information changes over time.
This model would surely coincide with DUE in the limit $\radius\to +\infty$, but, unfortunately, it seems to be quite complicated to be implemented in a numerical code.
For this reason, we have preferred to consider another scenario, described in the next section. 
We expect that, also in this case, convergence to DUE is guaranteed, but a precise theoretical framework for that is not easy to devise and it is out of the scope of this paper.

\subsubsection{The V2V-DUE model}\label{sec:V2V-DUE}
We assume that at each rendez-vous, each car $\c$ shares not only its current position and its destination, but also \textit{its planned path}, i.e.\ the function $(t,\j)\to\nextroad_\c(t,\j)$. 
Exploiting this information, each car can easily predict the path of the other cars it knows without arbitrariness (cf.\ Remark \ref{rem:discrepancy}), and then it can optimize its own path on the basis of that prediction.

For better clarity, let us consider the limit case $\radius\to +\infty$:  
in this case all cars re-encounter all other cars every $\delta c$ time units.
At every fixed time $\bar t$, each car $\c$ is able to forecast the positions of all other cars at any future time $t>\bar t$ using their (communicated) best strategies, and then it is able to devise its own best strategy for the current and future time.
At the next communication time, all the updated best strategies are shared with all the others in the new rendez-vous, then new predictions are made and finally new best strategies are computed.
This scenario actually corresponds to a \textit{delayed version} of the forward-backward iterations strategy used for computing DUE (Section \ref{sec:DUE}), where $\delta c$ time units are interposed between two consecutive forward-in-time dynamics prediction.
Roughly speaking, it is like cars start moving while they are still computing the equilibrium dynamics, because they have to wait $\delta c$ time units to get the others' up-to-date best strategies.

Summarizing, at any fixed time $\bar t$ and for any car $\c$:
\begin{enumerate}
        \item car $\c$ forgets all known cars encountered more than $\delta m$ time units in the past;
	\item positions of all cars $\c^\prime\in\mathcal K_\c(\bar t)$ are nowcast, in an imaginary fictitious simulation, from the last rendez-vous time $\hat t$ to time $\bar t$ using \eqref{ODE}-\eqref{def:vel} and $\nextroad_{\c^\prime}(s,\cdot)$, $s\in[\hat t,\bar t)$, assuming that only cars in $\mathcal K_\c(\bar t)$ are present on the network; 
	\item if $\delta c$ time units have passed since the last communications, car $\c$ checks for new rendez-vous. 
	It gets the current position, destination, and the function $\nextroad_{\c^\prime}(\cdot,\cdot)$ of any other car $\c^\prime$ within communication distance $\radius$. 
	If car $\c^\prime$ was already previously encountered, new pieces of information overwrite the old ones;
	\item positions of all cars $\c^\prime\in\mathcal K_\c(\bar t)$ known by car $\c$ are predicted at any time $t>\bar t$ using \eqref{ODE}-\eqref{def:vel} and $\nextroad_{\c^\prime}$, assuming that only cars in $\mathcal K_\c(\bar t)$ are present on the network;
	\item weights $w_\c(t,\cdot)$, $\forall t\geq \bar t$, are updated using cars in $\mathcal K_\c(\bar t)$ only;
	\item $\nextroad_\c(t,\cdot)$, $\forall t\geq \bar t$, is recomputed using \eqref{PPD-DUE} (in the time interval $[\bar t,\Tfin]$);
	\item $\nextroad_\c(\bar t,\cdot)$ is used in \eqref{ODE}-\eqref{def:vel} to move car $\c$ ahead.
\end{enumerate}

\begin{remark}\label{V2VDUE-to-DUE}
The model described above is expected to converge (in some sense) to the DUE for $\radius\to +\infty$ and $\delta c\to 0$.
Note also that in a time-discrete context, the time step $\Delta t$ plays the role of a delayer parameter as $\delta c$ does. 
In fact, since rendez-vous happen at each time step, the minimum delay is $\Delta t$ even if $\delta c=0$. 
Thus, considering time discretization too, also $\Delta t\to 0$ is obviously required for convergence.
\end{remark}

\subsubsection{Test 5}
\begin{figure}[ht!]
    \centering
    \subfloat[$3\times 3$, fixed OD pairs]{\label{subfig:TTT(R)-V2V-RUE-DUE-3x3-fixed}
    \includegraphics[width=.4\linewidth]{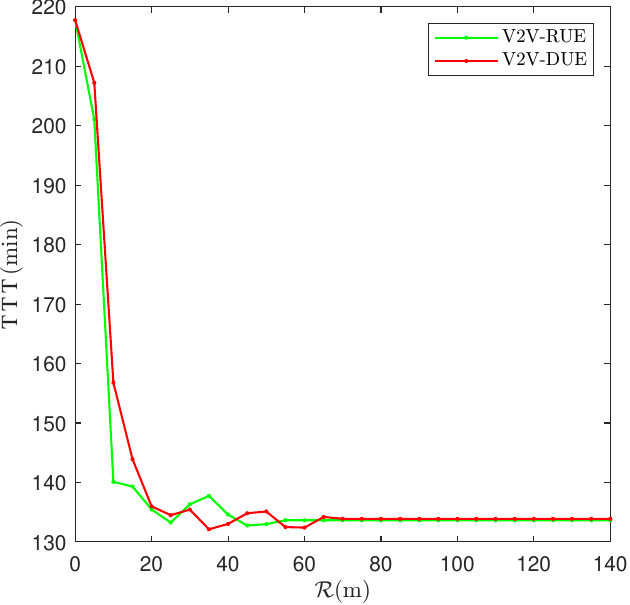}}\quad
    \subfloat[$3\times 3$, random OD pairs]{\label{subfig:TTT(R)-V2V-RUE-DUE-3x3-random}
    \includegraphics[width=.4\linewidth]{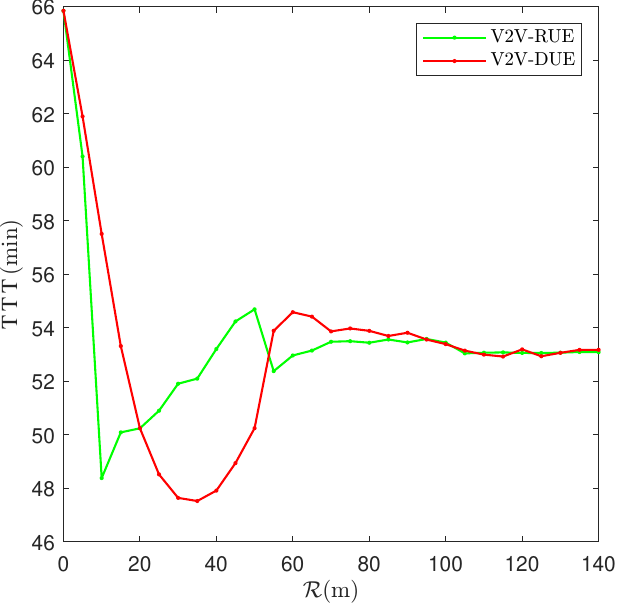}}\\
    \subfloat[$5\times 5$, random origins, fixed unique destination]{\label{subfig:TTT(R)-V2V-RUE-DUE-5x5-fixedrandom}
    \includegraphics[width=.4\linewidth]{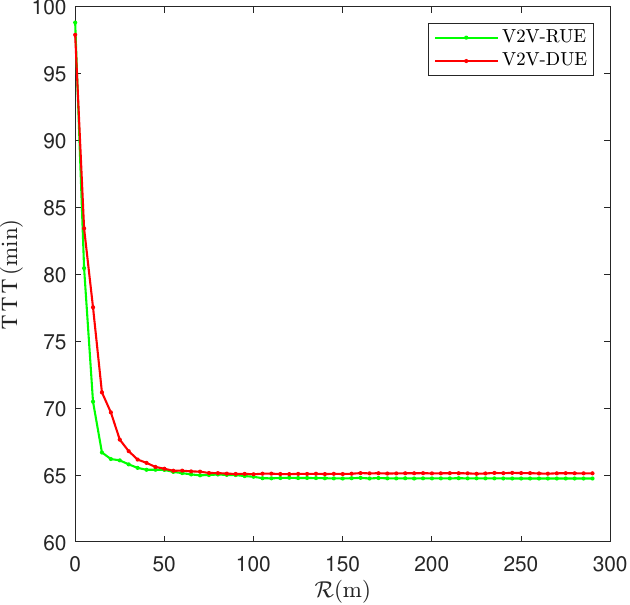}}\quad
    \subfloat[$5\times 5$, random OD pairs]{\label{subfig:TTT(R)-V2V-RUE-DUE-5x5-random}
    \includegraphics[width=.4\linewidth]{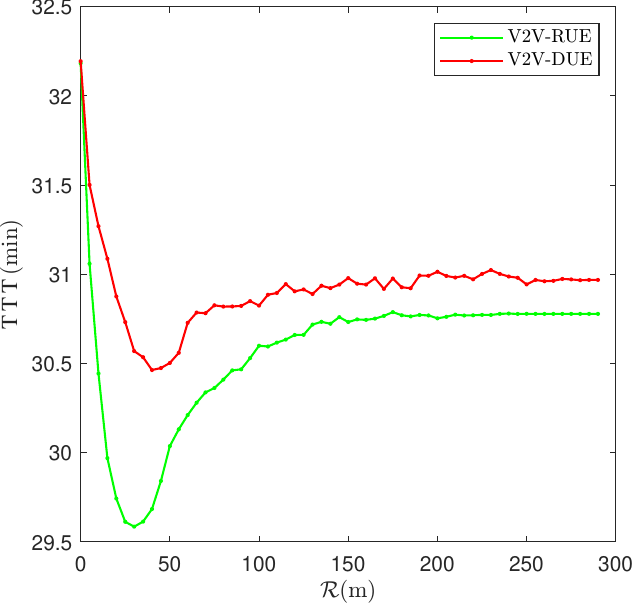}}\par
    \caption{Tests 3 and 5: $\TTT$ as a function of the communication range $\radius$, for V2V-RUE and V2V-DUE, Manhattan-like network $3\times 3$ and $5\times 5$, 100 cars}
    \label{fig:TTT(R)-V2V-RUE-DUE}
\end{figure}

In this section we repeat the same investigations as in Test 3 for V2V-DUE, with same parameters; we also compare V2V-RUE with V2V-DUE.
In Figs.\ \ref{subfig:TTT(R)-V2V-RUE-DUE-3x3-random}-\ref{subfig:TTT(R)-V2V-RUE-DUE-5x5-random} we observe again the blissful ignorance paradigm in the case of random OD pairs. This time, the capacity of forecasting the others' dynamics seems to mitigate the advantage of having a small communication range in the $5\times 5$ network, where traffic never reaches high degree of congestion, while it is even accentuated in the $3\times 3$ network where long queues appear.
In the case of fixed destinations or fixed OD pairs, instead, the $\TTT$ decreases almost monotonically as $\radius$ increases, see Fig.\ \ref{subfig:TTT(R)-V2V-RUE-DUE-3x3-fixed}.

We have also checked the correctness of the limit behaviors $\radius=0$ and $\radius=+\infty$. 
In the first case, the $\TTT$ of both V2V-RUE and V2V-DUE always coincides with that of BB.
In the second case, the $\TTT$ of V2V-RUE always coincides with that of RUE, while for V2V-DUE the situation is more delicate.
Whenever DUE exists and forward-backward iterations converge, the $\TTT$ of V2V-DUE coincides with that of DUE. 
If instead DUE is unstable and forward-backward iterations do not converge, we are not able to say anything about the convergence of V2V-DUE to DUE. 

Finally, let us note that, in some cases, the $\TTT$ of V2V-DUE is equal or larger than the $\TTT$ of \REV{V2V-RUE}, see Fig.\ \ref{subfig:TTT(R)-V2V-RUE-DUE-3x3-fixed}-\ref{subfig:TTT(R)-V2V-RUE-DUE-5x5-fixedrandom}.
Actually, we have also observed that in the setting of Fig.\ \ref{subfig:BB_vs_RUE_vs_DUE_random}, for some particular initial conditions, also DUE's $\TTT$ is equal or larger than RUE's $\TTT$. 
This is not in principle impossible to happen, and we think that the main reason is the great number of equivalently optimal paths in the Manhattan-like network. 


\section{Conclusions and future work}
In this paper we have proposed two models for traffic flow on network based on V2V communications. 
Unlike the majority of existing models, here vehicles share their destinations and possibly their planned path, thus allowing vehicles to predict the others' paths within a certain degree of accuracy. 
Since both V2V-RUE and V2V-DUE involve prediction of the others' behavior to a certain extent, in both cases we face a differential game which leads to a possible equilibrium configuration. 
This equilibrium is difficult to characterize and strongly depends on the system at hand, especially on the shape of the road network and OD pairs.

One of the most interesting results regards the role played by the communication range $\radius$. 
In the two extreme cases $\radius=0$, $\radius=+\infty$, we easily recover the well known behaviors BB and RUE/DUE, respectively. 
In between, where only limited information is available, new dynamics and traffic equilibria appear, which would be nice to be analytically characterized. 

\medskip

As future work, we aim at studying the many-particle limit of the V2V dynamics, in order to find the macroscopic counterpart of the proposed models.
This seems to be a quite challenging goal since classical multi-population models are not detailed enough to describe the dynamics. Even if we use separate density functions $\rho_\d(t,x)$ for describing groups of vehicles with different destinations, it remains the problem of keeping track of the different information collected by different vehicles during the journey.

\medskip 

\acknowl{
We would like to warmly thank Simone Cacace, Davide Vergni, and Alessio Olivierio for their useful comments and suggestions. 
We also want to send a very special thank to Maurizio Falcone, who introduced us to the beauty of the DPP.
}

\funding{
E.C.\ would like to thank the Italian Ministry of University and Research (MUR) to support this research with funds coming from PRIN Project 2022 PNRR entitled ``Heterogeneity on the road - Modeling, analysis, control'', No. 2022XJ9SX, and PRIN Project 2022 entitled ``Optimal control problems: analysis, approximation and applications'', No. 2022238YY5. \\
This study was carried out within the Spoke 7 of the MOST -- Sustainable Mobility National Research Center and received funding from the European Union Next-Generation EU (PIANO NAZIONALE DI RIPRESA E RESILIENZA (PNRR) – MISSIONE 4 COMPONENTE 2, INVESTIMENTO 1.4 – D.D. 1033 17/06/2022, CN00000023). 
This manuscript reflects only the authors' views and opinions. Neither the European Union nor the European Commission can be considered responsible for them. \\
E.C.\ and F.L.I.\ are funded by INdAM--GNCS Project, CUP E53C23001670001, entitled ``Numerical modeling and high-performance computing approaches for multiscale models of complex systems''.\\
F. L. I.\ is funded by INdAM--GNCS Project, CUP E53C24001950001, entitled ``Metodi avanzati per problemi di Mean Field Games ed applicazioni''.\\
F.L.I.\ is funded by Sapienza - University of Rome, project ``Partial differential equations towards control theory and climate modeling'' (project code: RM124190DEC62D0A).\\
E.C.\ and F.L.I.\ are members of the INdAM research group GNCS.
}


\complia{On behalf of all authors, the corresponding author states that there is no conflict of interest.}

\ethical{Not applicable.}


\appendix

\section{Implementation details}
The numerical codes for V2V-RUE and V2V-DUE are complex and tricky, and many details can be missed at first glance.
In this Appendix, we list the most important ones:
\begin{itemize}
    \item In the dynamical system \eqref{ODE}, the computation of $\delta_\c$ requires the computation of $\nextcar_\c$, which, in turn, requires the knowledge of $\nextroad_\c$. This is because the car ahead of $\c$ must be searched along the path of $\c$, therefore the path of $\c$ must be known;
    \item when method M4 is considered for DUE, the DPP must be slightly modified with respect to the standard version \eqref{PPD-DUE}. In fact, when one simulates the dynamics of the probe vehicle from the beginning to the end of any road $\r$, one must consider all the possible $\nextroad$'s it can have once it reaches the end of the road (because different $\nextroad$ means different $\nextcar$, then different dynamics, as pointed out in the previous comment). As a consequence, the function $w$ depends not only on $t$ and $\r$, but also on any road outgoing from the end of road $\r$. The minimum in \eqref{PPD-DUE} must be then extended by covering these additional cases;
    \item cars do not know themselves, i.e.\ $\c\notin \mathcal K_\c$. Moreover, any time a fictitious simulation is run (e.g., nowcast for V2V-RUE, nowcast and forecast for V2V-DUE), cars in $\mathcal K_\c$ cannot take into account the presence of $\c$ itself in their dynamics, otherwise inconsistencies can arise. 
    \item Once a car reaches its destination, it becomes inactive and it is no longer seen by the other, still active, cars. This is true both in real and fictitious simulations. Moreover, weights $w$ must be computed only considering active cars at any given time. Therefore, it can happen that in fictitious simulations, a car $\c^{\prime\prime}\in\mathcal K_\c \cap \mathcal K_{\c^\prime}$ can be, say, active for $\c$ and inactive for $\c^\prime$. Therefore, it is important to keep track of the status of each car at any time in all simulations in which it is involved.
\end{itemize}

\bibliographystyle{emibibliostyle}
\bibliography{biblio}
\end{document}